\documentclass[10pt]{amsart}
\usepackage[dvips]{epsfig}
\usepackage{amssymb}
\usepackage{ifthen}
\newlength{\widthofspace}
\newcommand{\negativespace}{\hspace{-\widthofspace}}

\newcommand{\skipthistext}[1]{}
\newcommand{\remark}[1]{}

\newcommand{\opname}[2][it]%
        {%
	\ifthenelse{\equal{#1}{it}}
	{
	\expandafter\def\csname #2\endcsname{\ensuremath{\mathit{#2}}}%
	}{}
	\ifthenelse{\equal{#1}{bb}}
	{
	\expandafter\def\csname #2\endcsname{\ensuremath{\mathbb{#2}}}%
	}{}
	\ifthenelse{\equal{#1}{rm}}
	{
	\expandafter\def\csname #2\endcsname{\ensuremath{\mathrm{#2}}}%
	}{}
	\ifthenelse{\equal{#1}{bf}}
	{
	\expandafter\def\csname #2\endcsname{\ensuremath{\mathbf{#2}}}%
	}{}
	\ifthenelse{\equal{#1}{tt}}
	{
	\expandafter\def\csname #2\endcsname{\ensuremath{\mathtt{#2}}}%
	}{}
	\ifthenelse{\equal{#1}{cal}}
	{
	\expandafter\def\csname #2\endcsname{\ensuremath{\mathcal{#2}}}%
	}{}
	\ifthenelse{\equal{#1}{bit}}
	{
	\expandafter\def\csname #2\endcsname{\ensuremath{\text{\boldmath $#2$}}}%
	}{}
	}
\opname{rel}
\opname{card}
\opname{Nd}
\opname[bit]{T}

\renewcommand{\int}[1]{\ensuremath{\stackrel{\circ}{#1}}}
\let\diffeo=\cong
\let\iso=\cong
\let\phi=\varphi

\def\proof{\relax}
\renewcommand{\proof}{\par\noindent{\em Proof:\ }\relax}
\newcommand{\epsymb}{$\square$}
\newcommand{\eproof}{ \rule{0mm}{1mm}\ \hfill \epsymb}
\newcommand{\bcs}{\,\natural\,}

\newcommand{\hcap}%
	{
  	\setlength{\unitlength}{0.01em}
  	\begin{picture}(100,58)
  	\put(50,0){\oval(50,60)[t]}
  	\end{picture}
	}
\newcommand{\boldcup}{\mbox{\boldmath $\cup$} }

\let\bd = \bolddelim

\newenvironment{theorem}[1][\negativespace]%
	{\vspace{5mm}\par\noindent {\bf Theorem #1. }\it}%
	{\rm\vspace{5mm}\par}

	{\vspace{5mm}\par\noindent {\bf Definition #1. }\it}%
	{\rm\vspace{5mm}\par}

\newenvironment{corollary}[1][\negativespace]%
	{\vspace{5mm}\par\noindent {\bf Corollary #1. }\it}%
	{\rm\vspace{5mm}\par}

	{\vspace{5mm}\par\noindent {\bf Lemma #1. }\it}%
	{\rm\vspace{5mm}\par}

	{\vspace{5mm}\par\noindent {\bf Question #1. }\it}%
	{\rm\vspace{5mm}\par}

	{\vspace{5mm}\par\noindent {\bf Remark #1. }\it}%
	{\rm\vspace{5mm}\par}

\newcommand{\drawpict}[1]%
	{\centerline{\epsfig{file=#1}}}

\newcommand{\putfigure}[3][h]%
	{
	\begin{figure}[#1]
	\drawpict{#2l}{}
	\caption{#3}\label{figure:#2}
	\end{figure}
	\expandafter\def\csname #2\endcsname{\ref{figure:#2}}
	\negativespace\negativespace
	}

\newcommand{\subtitle}[1]%
	{\vspace{5mm} \par\noindent {\bf\large {#1}\par}}

\let\Pn=1
\let\Pos=2
\opname{Wh}
\opname{tb}
\opname{rot}
\opname{tb}
\opname{id}
\opname[cal]{D}

\title{ A Convex decomposition theorem for four-manifolds}
\author{S. Akbulut and R. Matveyev}
\date{May 1997}

\begin{document}

\begin{abstract}
In this article we show that every smooth closed oriented
four-manifold admits a decomposition into two submanifolds
along common boundary. Each of these submanifolds is a
complex manifold with pseudo-convex boundary.
This imply, in particular, that every  smooth closed
simply-connected four-manifold is a Stein domain in the the
complement of a certain contractible 2-complex.
\end{abstract}

\maketitle

\subtitle{1. Introduction}
Exact manifold with pseudo-convex boundary
(PC manifold, for short)
is a compact complex manifold $X$, which admits strictly
pluri-subharmonic Morse function $\psi$, such that set
of maximum points of $\psi$ coincides with the boundary
$\partial X$. We prefer the term PC manifold, since
combination of words
``{\em compact Stein} manifold'' is likely to precipitate
heart palpitations in  some mathematicians.

Such manifold admits a symplectic structure
$\omega = \frac{i}{2} \partial\overline\partial\psi $
and serves as an analogue of closed symplectic
manifold.
Boundary $\partial X$ of PC manifold $X$ inherits
a contact structure $\xi$, which, in this case, 
is a distribution of
maximal complex subspaces in $\T X$ tangent to
$\partial X$.
In dimension four analogy between PC manifolds and
closed symplectic manifolds
is further illustrated by
the following two theorems in terms of
 Seiberg-Witten invariants.

\begin{theorem}[1a] {\em (C. Taubes, \cite{T})}
Let $(X,\omega)$ be a closed, symplectic four-manifold
and $K$ be Chern class of the canonical bundle
of an almost
complex structure compatible with $\omega$.
Then $SW_X( K ) = \pm 1$.
\end{theorem}

In the relative case Kronheimer and
Mrowka proved the following theorem:

\begin{theorem}[1b]
{\em (P. Kronheimer, T. Mrowka, \cite{KM})}
Let $(X,\omega)$ be compact, symplectic manifold and
$\xi$ be
a positive contact structure
on $\partial X$ compatible with $\omega$.
Then $SW_X( K ) = 1$, where $K$ is the canonical class
of $\omega$.
\end{theorem}

In particular, it was shown by P.~Kronheimer
and T.~Mrowka that properly embedded surface in PC
manifold satisfies Eliashberg-Ben\-ne\-quin inequality
analogous to adjunction inequality in the case of
closed symplectic manifold, once proper
conditions on the boundary of the surface are
imposed.
Namely,
if $F\subset X$ is a properly embedded surface,
such that
$\alpha=\partial F\subset\partial X$ is a Legendrian knot
with respect to induced contact structure on $\partial X$
and $f$ is
framing on $\alpha$ induced by a trivialization of
the normal bundle of $F$ in $X$, then
\[
[ \tb(\alpha)-f ]+|\rot(\alpha,F)|\leq -\chi(F),
\]
where $\tb(\alpha)$ is Thurston-Bennequin framing  defined
by the vector field along $\alpha$ transversal to
$\alpha$ and tangent to contact distribution.
To define rotation number observe that
$\T X|_{\partial X} \iso \xi \oplus \mathbb C$ and thus
$\xi \iso \Lambda^2 \T X|_{\partial X}$.
This allows us to view $\tb(\alpha)$ as a section of
$\Lambda^2 \T X|_\alpha$.
Then $\rot(\alpha,F)$ of Legendrian knot $\alpha$
bounding oriented proper surface $F \subset X$
is defined to be an obstruction to extend vector
field $\tb(\alpha)$ to a non-vanishing section
of $\Lambda^2 \T X |_F$, i.e.
\[
\rot(\alpha,F) =
c_1(\Lambda^2 \T X |_F,\tb(\alpha)) \hcap [F,\partial F].
\]
Details could be found, for example, in \cite{AM}.

\subtitle{2. Whitehead multiple of a knot}

Suppose $(K,f)$ is a framed knot in an oriented
3-manifold $M$.
{\em Positive Whitehead multiple} $P_n(K,f)$ of knot
$K$ is a band connected
sum of $n$ parallel (according to framing $f$)  copies
of $K$ equipped
with alternating orientations as on Figure \Pn\ 
 (we assume usual orientation of $\mathbb R^3$).
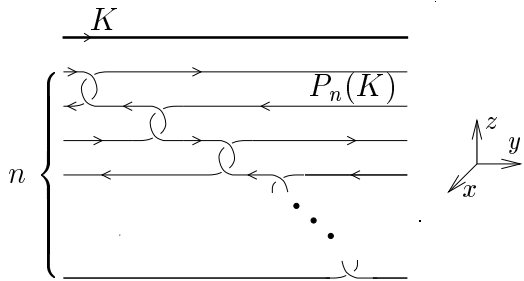
\putfigure{Pn}{Positive Whitehead double of knot $K$}
We shall omit framing $f$ from the notation when it  is
clear from the
context or irrelevant to the discussion.

This construction is a generalization of Whitehead
double of a
knot, namely $P_2(K)=\Wh(K)$.

Here we list some properties of the Whitehead multiple
$P_n(K,f)$.

{\bf 1.}
The homology class of $P_n(K,f)$ is either that of $K$
if $n$ is odd, or zero if $n$ is even.
In fact, for odd values of $n$, $P_n(K,f)$ is homotopic
to $K$.
In both cases,
$P_n(K,f)$ can be naturally equipped with a framing,
which we also call $f$. It is framing corresponding to
the given framing of $K$
in case of odd $n$ (this is defined correctly because $K$
and
$P_n(K,f)$ are in the same homology class). Canonical
framing
of $P_n(K,f)$ for even $n$ is just a zero framing, which
is
well-defined for the closed curve homologous to zero.

{\bf 2.} The key property of the above operation is
that given a Legendrian knot $K$ in tight contact
3-manifold
with zero Thurston-Bennequin invariant
\skipthistext{$\tb(P_n(K,0))=n-1$.}
we can produce another Legendrian knot in small
neighborhood of $K$ homotopic to $K$ and with an
arbitrary large
Thurston-Bennequin invariant.
More precisely, let $(K,f)$  be framed Legendrian knot in
tight contact manifold $M$.
Let $K_1, \ldots , K_n$ be parallel copies of $K$
corresponding to framing $f$.
Isotopy class of link $(K_1, \ldots , K_n)$
has a Legendrian representative
with $K_1=K$ and $\tb(K_i)-f = -|\tb(K)-f|$, $i=2,\ldots, n$.
This is obvious for $K$ being Legendrian unknot with
$\tb(K) = 1$ and $\rot(K) = 0$.
The general case follows from the fact that any two Legendrian
knots have contactomorphic neighborhoods.
Different proof is given in \cite{AM1}.

Each band in the construction of $P_n(K, \tb(K))$ contributes
$1$ into $\tb(P_n(K, \tb(K)))$, thus we have
\[
\tb \bd( P_n(K, \tb(K)) \bd) - \tb(K) = n-1.
\]

One has to
choose $n$
to be odd to make $P_n(K,f)$ homotopic to $K$.

{\bf 3.}
Suppose knot $K\subset M^3$ is a boundary of properly embedded
disc $D$ in
smooth 4-manifold $X$ with
$\partial X = M^3$.
Taking a ribbon sum of $n$ copies of disc $D$,
using the same ribbons, which are used to construct
$P_n(K)$,
we obtain
another disc $P_n(D)$.
It is bounded by $P_n(K,f)$, where $f$ is the framing
induced by
trivialization of the normal bundle $\nu_X (D)$ of $D$ in
$X$.
Note also, that disc $P_n(D)$ can be constructed in
arbitrary
small neighborhood of $D$.

{\bf 4.}
If $(l_1,l_2)$ is a Hopf link then
link $( P_n(l_1,0) , l_2 )$ in $S^3$ is symmetric
(i.e. there is a diffeomorphism of $S^3$ interchanging
components of the link).
Moreover, two triads  $( B^4 , P_n(D_1) , D_2 )$ and
$( B^4 , D_1, P_n(D_2) )$ are diffeomorphic.
Here, $D_1$ and $D_2$ is a pair of linear 2-discs in
$B^4$ intersecting at one point.

Paying homage to the popularity of  physics terminology
first author suggested to call manifold $W_n$ on
Figure \Pos\
a positron.
Manifold $W_n$ is obtained
from $B^4$ by removing open regular neighborhood of
$D_1$ (this has an effect of turning $B^4$ to
$S^1\times B^3$)
and then attaching a 2-handle to the framed knot
$P_n(\partial D_2)$.
It is contractible if $n$ is odd and has PC
structure if $n\ge 2$.
\putfigure{Wn}{Positron}

\subtitle{3. Handlebodies of PC manifolds}

Handlebodies of four-dimensional PC manifolds
are characterized by the
following theorem of Ya. Eliashberg.

\begin{theorem}[2]
{\em (Ya. Eliashberg, \cite{E}; see also \cite{G}) }
Let
$X=B^4\cup(\mbox{1-handles})\cup(\mbox{2-handles})$
be four-dimensional handlebody with one
0-handle and no 3- or 4-handles.
Then:
\begin{itemize}
\item The standard { PC} structure on $B^4$ can
      be extended over 1-handles, so that
      manifold $X_1=B^4\cup(\mbox{1-handles})$
      has pseudo-convex boundary.
\item If each 2-handle is attached to
      $\partial X_1$ along a Legendrian
      knot  with framing one
      less then Thurston-Bennequin framing
      of this knot, then the complex structure
      on $X_1$ can be extended over
      2-handles to a complex structure on $X$,
      which makes $X$ a { PC} manifold.
\end{itemize}
\end{theorem}

In this section we shall study ``partial handlebodies''
obtained by attaching 2-handles on top of a PC
manifold.
Let $Z$ be a PC manifold and $h$ be a two handle
attached to $\partial Z$ along a Legendrian knot $K\subset\partial Z$
with  framing $f$.
If $\tb(K) \ge f+1$ then by $C^0$-small smooth isotopy of $K$
we can decrease Thurston-Bennequin invariant of $K$ and
make it equal to $f+1$.
Therefore, by a theorem of Eliashberg,  manifold
$Z \cup h$ possesses PC structure.

However, in general, it is not possible to increase
Thurston-Bennequin invariant of $K$ by isotopy,
and equip $Z \cup h$ with
the structure of PC  manifold.
Thus, we make the following definition: the {\em defect} $\D(h^2)$ of a
2-handle $h^2$ attached to a Legendrian knot $K$ on the
boundary of PC  manifold with framing $f$ is a number
$\max\{f+1-\tb(K),0\}$.
If we have several 2-handles $h_1^2, \ldots, h_n^2$ attached to a Legendrian
link on the boundary of PC  manifold $Z$,
then the defect $\D(Z \cup \boldcup_i h^2_i)$ of this partial handlebody
built on
top of $Z$ is a sum
of the defects of individual handles.
So if the defect is zero, the PC structure extends over 2-handles.
Given such partial handlebody
$Z \cup \boldcup_i h_i$ with a base $Z$ being a PC manifold ,
we can build another $Z' \cup \boldcup_i h'_i$ with lesser
defect, so that $Z$ and $Z \cup \boldcup_i h_i$ are homotopy
equivalent to $Z'$ and $Z' \cup \boldcup_i h'_i$,
respectively.
\putfigure{Zh}{Manifold $Z\cup h_i$}
To see this, let $h_i$ be a handle with a non-zero defect.
Let $D_i$ be a cocore of handle $h_i$ and $m_i$ be a
meridian of $P_n(\partial D_i)$ in $\partial Z$, see Figure \Zh.
We may view $D_i$ as a disk in $Z$.
Note that $(P_n(\partial D_i), m_i)$ is a small
(in a chart) Hopf link in $\partial Z$.
Manifold $Z'$ is obtained by removing
$\Nd_Z(P_n(D_i))$ ---
a tubular neighborhood of $P_n(D_i)$ in $Z$,
from $Z$ and attaching a
2-handle to $P_k(m_i)$.
We assume that $k$ is odd and $k\ge 3$.
Manifold $Z'$  is boundary connected sum of
$Z$ and a positron $W_k$, hence it is PC manifold
homotopy equivalent to $Z$.
New handle $h'_i$ is attached to the connected
sum of attaching circle of $h_i$ and
$P_n(\text{core of 1-handle in positron})$,
thus it's defect is $n$ less than defect of
$h_i$.
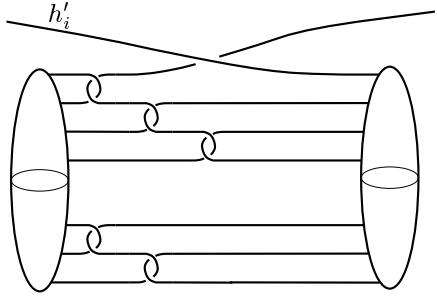
\putfigure{Zhprime}{Manifold $Z'\cup h'_i$ for $n=4$, $k=3$.}
Manifold $Z'\cup h'_i$ is shown on Figure \Zhprime.

\subtitle{4. Convex Decomposition Theorem}

We will use above construction to prove the following
convex decomposition theorem.

\begin{theorem}[3] Let $X = X_1 \cup_\partial X_2$ be a
decomposition of a closed smooth oriented 4-manifold
into a union of two compact, smooth, codimension zero
submanifolds $X_1$
and $X_2$ along common boundary.
Suppose each $X_i$, $i=1,2$, has a handlebody
without 3- and 4-handles.
Then there exist another decomposition
$X = \tilde{X}_1 \cup_\partial \tilde{X}_2$,
such that manifolds
$\tilde{X}_1$ and  $-\tilde{X}_2$ admit structures of PC
 manifolds and each $\tilde{X}_i$ is homotopy equivalent to
$X_i$, $i=1,2$.
\end{theorem}

\noindent{\em Proof:}
Consider handlebodies of $X_i$, $i=1,2$, with the
properties stated in the assumption of the theorem.
Let $Y_i$ be a union of 0- and 1-handles in $X_i$.
According to theorem of Eliashberg (Theorem 2, above) $Y_i$
is a PC  manifold.
Complex structures on $Y_i$ are chosen so that
complex orientation coincides with the induced
orientation on $Y_1 \subset X$ and is opposite on
$Y_2 \subset X$.
Since every curve in contact manifold is isotopic to
a Legendrian curve via smooth $C^0$-small isotopy,
we can assume that 2-handles in handlebodies of
$X_i$ are attached to Legendrian knots in
$\partial Y_i$.
Let $h$ be a 2-handle in $X_1$ with a non-zero
defect, $D$ be the cocore of $h$, $d = \partial D$,
 $m$ be the meridian
of $P_n(d)$ and $F$ be  a trivial embedded disc in $X_2$
bounded by $m$.

As in the construction in previous section,
we remove $P_n(D)$
from $X_1$, and attach handle to $P_k(m)$
with framing $0$,
reducing defect of $h$ and, therefore, total defect of
$X_1$ by $n$.
Manifold $X'_1$ can be built inside of $X$, namely
\[
X'_1 = [ \, X_1 \setminus \Nd_{X_1} ( P_n D ) \, ] \cup
\overline{\Nd_{X_2} ( P_k  F )}.
\]
Here $\Nd_X (Y)$ stands for a tubular neighborhood
of $Y$ in $X$.
Its complement $X'_2$ is obtained from $X_2$ by
attaching a new 2-handle $g = \overline{\Nd_{X_1}( P_n D )}$ and removing
neighborhood of $P_k(F)$.
Since $X_1$ and $-X_2$ induce the same orientation on
their common boundary, positive Whitehead multiple
is the same whether it is considered in
$\partial X_1$ or $\partial (-X_2)$.
If we choose $n$ to be defect of $h$ and $k$ to
be odd and greater then or equal to defect of $g$
(after Legendrianization of attaching circle),
than total defect of $X_1$ is reduced by $n$ and
defect of $X_2$ is not increased.
By applying this procedure to every 2-handle of $X_1$
with non-zero defect, we obtain
manifold $\bar{X}_1$ with pseudo-convex boundary,
and it's complement $\bar{X}_2$
has a defect less then or  equal to the defect of $X_2$.
To finish the
proof one has to apply the same procedure to decomposition
$-X=(-\bar{X}_2) \cup (-\bar{X}_1)$ to obtain decomposition
$X=\tilde{X}_1 \cup_\partial \tilde{X}_2$,
with $\tilde{X}_1$ and $-\tilde{X}_2$ being
PC manifolds homotopy equivalent to
$X_1$ and $-X_2$, respectively.\eproof

\begin{corollary}
Every closed simply-connected four-manifold $X$ possesses a
structure of complex manifold with pseudo-convex boundary
 in the complement of some compact contractible
submanifold (which is also PC manifold).
\end{corollary}
\proof
Consider arbitrary handle decomposition of $X$.
Let $Y_1$ be a union of all 0- and 1-handles in $X$.
The fundamental group of $Y_1$ is free and the natural map
$\pi_1(\partial Y_1)\rightarrow \pi_1(Y_1)$  is an isomorphism.

We shall rearrange handlebody of $X$ by introducing pairs of dual 2-
and 3-handles
and handle addition, so that
among 2-handles of new handlebody of $X$ it is possible to  choose
 a set $\{h_i\}$ such that
$X_1=Y_1\cup\boldcup_i h_i$ is a contractible manifold.

More detailed description of this construction is as follows:
Let $\{x_1,\ldots ,x_l\}$ be a free basis of $\pi_1(\partial Y_1)$.
If we fix paths from a base-point to attaching spheres of 2-handles then
they represent elements of $\pi_1(\partial Y_1)$, say $y_1,\ldots ,y_L$.
Since $X$ is simply-connected, $\{y_1,\ldots ,y_L\}$ normally generate
$\pi_1(\partial Y_1) \diffeo \pi_1(Y_1)$.
Thus, each $\tilde{x}_i$, $i=1, \ldots, l$
is a product of elements adjoint to
$y_1,\ldots ,y_L$.
Introduce a pair of canceling 2- and 3-handles.
We can slide this new 2-handle over handles
corresponding
to the elements $y_1,\ldots ,y_L$ in the
decomposition of $\tilde{x}_i$,
so that its attaching circle become homotopic
to $\tilde{x}_i$.
Union $X_1$ of manifold $Y_1$ and 2-handles
obtained by construction above for every element
in $\{\tilde{x}_1,\ldots ,\tilde{x}_l\}$ is a
contractible submanifold of $X$.

Manifolds $X_1$ and
$X_2=\overline{X\setminus X_1}$
have only handles of indexes
less then 3 in their handle decompositions.
Thus, Theorem 3 can be applied to decomposition
$X=X_1 \cup_\partial X_2$, which finishes
proof of the corollary.\eproof

Above corollary implies that every
closed simply-connected manifold
posseses structure of Stein domain in the
complement of certain contractible two-dimensional complex,
hence every closed embedded surface $F$ in
the complement of this
complex satisfies adjunction inequality:
\[
-\chi (F) \geq F \cdot F + K \hcap F.
\]

\subtitle{5. Corks with pseudo-convex boundary}

\begin{theorem}[4]
{\rm(See \cite{CH, M})}
Let $M_1$, $M_2$ be two closed, smooth, simply-connected,
h-cobordant 4-manifolds.
Then:
\[
M_1=N\cup_{\phi_1} A_1 , \hspace{5mm}
M_2=N\cup_{\phi_2} A_2 ,
\]
where $N$ is simply-connected, $A_1$ and $A_2$ are contractible
and diffeomorphic to each other and
$\phi_i : \partial A_i \rightarrow \partial N$, $i=1,2$, are
some diffeomorphisms.
\end{theorem}

Manifolds $A_1$ and $A_2$ have come to be known as
{\em corks}, \cite{K}.

The proof of Theorem 3 can be adapted to show that
the corks and manifold $N$
in the theorem above can be made pseudo-convex.
\begin{theorem}[5]
Decompositions in Theorem 4 can be made pseudo-convex.
\end{theorem}

\proof
Start with decompositions provided by Theorem 4.
Note that by the construction, see \cite{M}, manifolds
$N$, $A_1$ and $A_2$ have handlebodies without handles of
indexes 3 and 4.
We divide the rest of the proof into three steps:

\vspace{0.5 em}\par
\noindent{\bf 1.} Make $N$ pseudo-convex.

Apply procedure from the proof of Theorem 3 to both
decompositions
$M_1 = N \cup_{\phi_1} A_1$, and
$M_2 = N \cup_{\phi_2} A_2$,
simultaneously.
We obtain new decompositions
$M_1 = N' \cup_{\phi'_1} A'_1$, and
$M_2 = N' \cup_{\phi'_2} A'_2$, where
$N'$ is pseudo-convex and homotopy equivalent to
$N$; $A'_1$, $A'_2$ are contractible, but not
necessarily diffeomorphic to each other.

\vspace{0.5 em}\par
\noindent{\bf 2.} Make $A'_1$ and $A'_2$ pseudo-convex.

Consider decompositions resulting from step 1 ``up side down'':
\begin{eqnarray*}
-M_1 & = & (-A'_1) \cup_{\phi'_1{}^{-1}} (-N') \text{  and}\\
-M_2 & = & (-A'_2) \cup_{\phi'_2{}^{-1}} (-N').
\end{eqnarray*}

Suppose $h$ is a $2$-handle in handlebody of $-A'_1$ with
non-zero defect.
Let $D$ be a cocore of $h$,
$d = \partial D$,
$m$ be a meridian of $P_n(d)$
and $F$ be a trivial embedded disk in
$N$ bounded by $\phi'_1 (m)$
(note that $m$ is unknot in
$\partial A_1 \diffeo \partial N$).
We set
\begin{eqnarray*}
A''_1 & = & [A'_1 \setminus \Nd(P_n D)] \cup
\overline{\Nd(P_k F)},\\
N''   & = & [N' \setminus \Nd(P_k F)] \cup
\overline{\Nd(P_n D)}.
\end{eqnarray*}

Hence $M_1= N'' \cup A_{1}''$. Now we have to find manifold $A''_2$,
so that $ N''\cup A''_{2} = M_2$
and defect of $A''_2$ is not greater than
defect of $A'_2$.
Note that $N''$ is obtained from $N'$
by attaching a 2-handle along
$\phi'_1(P_n (d))$ and then removing
Whitehead multiple of its cocore.
Thus, $\partial N''$ is the result of
the surgery of $\partial N'$ along the link
$(l_1, l_2)$, where
$l_1 = \phi'_1( P_n(d) )$ and
$l_2 = P_k(\text{meridian of $l_1$})$.
Take
$$A''_2 = [A'_2 \cup h'] \setminus
\Nd( \, P_k(\text{cocore of $h'$}) \, ) , $$
where $h'$ is a 2-handle attached along
$\phi'_2{}^{-1} \circ \phi'_1 (P_n(d))$.
Boundary of $A''_2$ is obtained by the
surgery of $\partial A'_2$
along link
$(\phi'_2{}^{-1} (l_1),\phi'_2{}^{-1} (l_1))$,
therefore $\phi'_2$
extends to
$\phi''_2 : \partial A''_2 \rightarrow \partial N''$.
Manifold $N'' \cup_{\phi''_2} A''_2$ is
diffeomorphic to connected sum of
$N' \cup_{\phi'_2} A'_2$ and the double of positron $W_k$.
It is easy to see that the double of a
positron is diffeomorphic to $S^4$,
which implies that
$N'' \cup_{\phi''_2} A''_2 \diffeo N_2$.
Now, if we choose $k$ greater then
$\max\{ \D(\overline{\Nd(P_k F)}) , \D(h') \} $\ \
(here $\overline{\Nd(P_k F)}$ is considered as a 2-handles
attached to $N'$) and take $n$ to be the defect of $h$,
then defect of $A'_1$ is reduced by $n$
and defects of $A'_2$ and $N'_2$ are not increased.
Apply the above construction to every 2-handle
with non-zero defect in $A'_1$ and $A'_2$.
The resulting decompositions
\begin{eqnarray*}
M_1 & = & N'' \cup_{\phi''_1} A''_1 ,\\
M_2 & = & N'' \cup_{\phi''_2} A''_2
\end{eqnarray*}
are pseudo-convex, but corks are not diffeomorphic
to each other.
To fix that is the subject of the next step.

\vspace{0.5 em}\par
\noindent{\bf 3.} Make $A''_1$ diffeomorphic to $A''_2$.

It is shown in \cite{M} that manifolds
$A_1$, $A_2$ have the
property that their doubles
$A_i \cup_\partial (-A_i)$ as well as their union
$A_1 \cup_{\phi_1^{-1}\circ\phi_2} (-A_2)$
are diffeomorphic to $S^4$.
Not difficult but rather technical calculation shows
that this properties are preserved under
modifications from steps 1 and 2 above.
Thus we have
\begin{eqnarray*}
&&A''_i \cup_\partial(-A''_i) \diffeo S^4,\;i=1,2;\\
&&A''_1 \cup_{\phi_1^{-1}\circ\phi_2} (-A''_2) \diffeo S^4.
\end{eqnarray*}

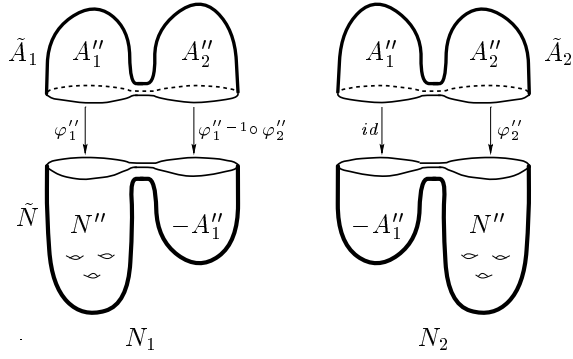
\putfigure{trick}{Convex decompositions with diffeomorphic corks.}

We define
\begin{eqnarray*}
\tilde N   & = & N''   \bcs (-A''_1), \\
\tilde A_1 & = & A''_1 \bcs A''_2,    \\
\tilde A_2 & = & A''_2 \bcs A''_1 ,
\end{eqnarray*}
where $X\bcs Y$ stands for boundary connected
sum of $X$ and $Y$.
Since boundary connected sum of PC manifolds
is a PC manifold,
$\tilde N$, $\tilde A_1$ and $\tilde A_2$ have
pseudo-convex boundary
and, obviously, $\tilde A_1 \diffeo \tilde A_2$.
Now we calculate
\begin{eqnarray*}
\tilde N \cup \tilde A_1
& = &
[N''   \bcs (-A''_1)]
\cup_{\phi_1 \bcs [\phi_1^{-1}\circ\phi_2]}
[A''_1 \bcs A''_2]\\
& \diffeo &
[N'' \cup_{\phi_1} A_1] \,\#\, [(-A''_1)
\cup_{\phi_1^{-1}\circ\phi_2} A''_2]\\
& \diffeo &
N_1 \,\#\, S^4 \diffeo N_1.
\end{eqnarray*}

Analogously,
\begin{eqnarray*}
\tilde N \cup \tilde A_2
& = &
[N''   \bcs (-A''_1)]
\cup_{\phi_2 \bcs \id}
[A''_2 \bcs A''_1]\\
& \diffeo &
[N'' \cup_{\phi_2} A_2] \,\#\, [(-A''_1)
\cup_{\id} A''_1]\\
& \diffeo &
N_2 \,\#\, S^4 \diffeo N_2.
\end{eqnarray*}
This gives us convex decomposition of $N_1$
and $N_2$ with diffeomorphic   corks
$\tilde A_1$ and $\tilde A_2$ and finishes
the proof of Theorem 5.\eproof

\subtitle{6. Questions and remarks}
We would like to conclude with some questions and remarks.

{\bf Question 1.}
Is it always possible to find convex decomposition as in Theorem 3,
so that
 contact structures on $\partial \tilde X_1$ and $\partial (-\tilde X_2)$
coincide.
(Authors can show that contact structures on
$\partial \tilde X_1$ and $\partial (-\tilde X_2)$
can be made homotopic.)

{\bf Question 2.}
Does Theorem 3 (and possibly positive answer to Question 1) pose
any restriction to the genera of embedded surfaces in four-manifold.


\begin {thebibliography}{MM1}

\bibitem[AM]{AM}	S. Akbulut, R. Matveyev,
			{\em  A note on Contact Structures.},
			Pacific J. of Math.

\bibitem[AM1]{AM1}	S. Akbulut, R. Matveyev,
			{\em Exotic structures and adjunction inequality.},
			Turkish J. of Math., vol. 21, No 1 (1997) pp. 47--52.

\bibitem[CH]{CH}	C.L. Curtis,  M.H. Freedman, W.C. Hsiang, R. Stong,
                    	{\em A decomposition theorem for h-cobordant smooth
                    	simply-connected compact 4-manifolds.},
		    	Invent. Math. {\bf 123} (1996) no. 2, pp. 343--348.

\bibitem[E]{E}      	Ya. Eliashberg,
			{\em Topological characterization of Stain
                    	manifolds of dimension $> 2$.},
			International J. of Math.
			Vol. 1, No 1 (1990) pp. 29-46.

\bibitem[G]{G}		R. Gompf,
			{\em Handlebody construction of Stein surfaces.},
			(preprint).

\bibitem[K]{K} 		R. Kirby,
			{\em Akbulut's corks and h-cobordisms
			of smooth, simply connected $4$-manifolds.},
			(preprint).

\bibitem[KM]{KM}	P. Kronheimer, T. Mrowka,
			{\em Monopoles and contact structures},
			preprint.

\bibitem[M]{M}		R. Matveyev,
			{\em A decomposition of smooth simply-connected
			$h$-cobordant $4$-manifolds.},
			J. Differential Geom. {\bf 44} (1996), no. 3, pp.
571--582.

\bibitem[T]{T}		C.H. Taubes,
			{\em The Seiberg-Witten invariants and
			symplectic forms.}
			Math. Res. Letters {\bf 1} (1994), pp. 808--822.

\end {thebibliography}
\end{document}